
\mag=\magstep1
\documentstyle{amsppt}

\topmatter
\title Subadjunction of log canonical divisors for a subvariety of
codimension 2
\endtitle
\author Yujiro Kawamata
\endauthor

\rightheadtext{adjunction}

\address Department of Mathematical Sciences, University of Tokyo, Komaba,
Meguro, Tokyo, 153, Japan \endaddress
\email kawamata\@tansei.cc.u-tokyo.ac.jp\endemail

\keywords adjunction, log canonical divisor, moduli of curves
\endkeywords

\abstract
We obtain a formula which relates the log canonical divisor of
the ambient space with that of a subvariety of codimension 2
by using Knudsen's moduli space of pointed stable curves of genus 0.
\endabstract

\endtopmatter

\document

\head Introduction
\endhead

The adjunction formula relates canonical divisors of varieties.
Let $X$ be a smooth variety and $S$ a smooth divisor.
Then we have $(K_X + S) \vert _S \sim K_S$.
But if $X$ has singularities, then we need to consider the
{\it subadjunction} as observed by M. Reid ([KMM, 5.1.9]).
For example, if $X$ is a singular conic in $\Bbb P^3$ and $S$
one of its generators, then $(K_X + S) \vert_S \sim_{\Bbb Q} - \frac 32 H$,
while $K_S \sim -2H$, where $\sim_{\Bbb Q}$ stands for the $\Bbb Q$-linear
equivalence and $H$ is the hyperplane section.
So it is more natural to consider the pair $(X, D)$ of the variety and
a $\Bbb Q$-divisor, where $D = S + \text{other components}$
with $S$ a prime divisor, and compare
$K_X + D$ and $K_S + D_S$ for some $\Bbb Q$-divisor $D_S$ on $S$.
We are also interested in the relationship between singularities
of the pair $(X, D)$ and those of the restricted pair $(S, D_S)$.
By using the residue map, we obtain a subadjunction theorem
(cf. [K3, Proposition 1.7]):
$(K_X + D) \vert_S \sim_{\Bbb Q} K_S + D_S$
for a canonically determined
effective $\Bbb Q$-divisor $D_S$ on $S$, and
if $(X, D)$ is LC (log canonical), then so is $(S, D_S)$.

If there are two smooth divisors $S_1, S_2$ on a smooth variety $X$
with a transversal intersection
$W = S_1 \cap S_2$, then by using the residue map twice,
we obtain $(K_X + S_1 + S_2) \vert_W \sim K_W$.
The purpose of this paper is to extend the subadjunction formula for a
codimension 2 subvariety $W$ even if there are no intermediate
divisors $S_j$, and answer [K3, Question 1.8] in the case of
codimension 2 (we use the notation of [KMM] and [K3]):

\proclaim{Theorem 1}
Let $X$ be a normal vairiety which has only KLT
singularities,
$D$ an effective $\Bbb Q$-Cartier divisor such that
$(X, D)$ is LC, and $W$ a minimal element of
$CLC(X, D)$ (the set of centers of log canonical singularities).
Assume that $\text{codim }W = 2$.
Then there exist canonically determined
effective $\Bbb Q$-divisors $M_W$ and $D_W$ on $W$
such that $(K_X + D) \vert_W \sim_{\Bbb Q} K_W + M_W + D_W$.
Moreover, if $D= D' + D''$ with $D'$ (resp. $D''$) the sum for irreducible
components which contain (resp. do not contain) $W$, then
$M_W$ is determined only by the pair $(X, D')$.
If $X$ is affine, then there exists an
effective $\Bbb Q$-divisor $M_W'$ such that $M_W' \sim_{\Bbb Q} M_W$ and
the pair $(W, M'_W + D_W)$ is KLT.
\endproclaim

The divisor $D_W$ is the {\it local contribution} and
appears by the same reason as in the case of divisorial subadjunction formula,
while the divisor $M_W$ is the {\it global contribution} which comes from
the moduli space of curves.
The latter is a new ingredient in the higher codimensional subadjunction.
We shall use Knudsen's moduli space of pointed stable curves of genus 0
in order to define $M_W$.
The naturally determined $\Bbb Q$-divisor
$M_W$ does not give us the KLT singularities, so we
need to introduce $M'_W$ which is not canonically determined.
We encountered a similar situation already in [N].
When $W$ an element of $CLC(X, D)$ which is not minimal, we still expect that
the pair $(W, M'_W + D_W)$ is LC for some $M'_W$.  In order to prove this,
one should show that $\Cal L$ in Theorem 4 is semi-ample.

Kodaira's canonical bundle formula for an elliptic surface $f: X \to B$
$$
K_X = f^*(K_B + \sum_i \frac{a_i}{12}P_i + \sum_j \frac {m_j-1}{m_j}Q_j +
\sum_k \frac {b_k}{12}R_k)
$$
can be regarded as an adjunction formula, where the $f^{-1}(P_i)$ are
the singular fibers of the types $II, III, IV, I^*, II^*, III^*$ or $IV^*$, and
the $f^{-1}(Q_j)$ and the $f^{-1}(R_k)$ are of the type $m_jI_{b_k}$.
The last sum $\sum_k \frac {b_k}{12}R_k$ is the global contribution in this
case. The main lemma (Theorem 4) is a log version of this kind of
formula which is obtained as the
{\it integration along fibers} (cf. Example 8).

\head 1. Proof of the main theorem
\endhead

We start with a lemma which gives a standard section of some multiple of
log canonical divisors.

\proclaim{Lemma 2}
Let $D = \sum_{j=1}^k d_jP_j$ be an effective $\Bbb Q$-divisor on
$C \cong \Bbb P^1$ such that
$0 < d_j \le 1$ for all $j$ and $\sum_j d_j = 2$.
Let $m$ be a positive integer such that $md_j \in \Bbb N$ for all $j$.
Then there exist a positive integer $p$ and a canonically defined section
$\omega_{D, m} \in \Gamma(C, mp(K_C + D))$
which depend only on $D$ and $m$.
\endproclaim

\demo{Proof}
We consider the set
$$
\align
P(D, m) = \{&J = (j_1, j'_1, \ldots, j_m, j'_m) \in \Bbb N^{2m} ; \\
&1 \le j_i, j'_i \le k, j_i \ne j'_i \text{ for all } i,
\#\{i \vert j_i = j \text{ or } j'_i = j\} = md_j\}.
\endalign
$$
Let $\omega_{j,j'}$ be a rational 1-form on $C$ such that
$\text{div}(\omega_{j,j'}) = P_j + P_{j'}$,
$\text{res}_{P_j}(\omega_{j,j'}) = 1$ and $\text{res}_{P_{j'}}(\omega_{j,j'})
= -1$.
Then we define
$$
\omega_J = \prod_{i=1}^m \omega_{j_i,j'_i}
\in \Gamma(C, m(K_C + D))
$$
and
$$
\omega_{D, m} = \prod_{J \in P(D,m)} \omega_J \in \Gamma(C, mp(K_C + D))
$$
where $p = \# P(D, m)$.
\hfill $\square$
\enddemo

Let $\Cal M_{0,n}$ be the moduli space of $n$-pointed stable curves of genus
$0$, $f_{0,n}: \Cal Z_{0,n} \to \Cal M_{0,n}$ the universal family,
and $\Cal P_1, \ldots, \Cal P_n$
the sections of $f_{0,n}$ which correspond to the
marked points ([Kn]).  $\Cal M_{0,n}$ is smooth.
Let $\Cal M_{0,n}^o$ be the open part which parametrizes the smooth curves.
$\Delta_{0,n} = \Cal M_{0,n} \setminus \Cal M_{0,n}^o$
is a normal crossing divisor
$\Delta_{0,n} = \sum_S \Delta_S$, where the $S = \{S',S''\}$
run all the decompositions
$\{1, \ldots, n\} = S' \cup S''$ such that $\# S', \#S'' \ge 2$.
$f_{0,n}^*\Delta_S$ has 2 irreducible components $\Cal F_{S'}, \Cal F_{S''}$,
to which the $\Cal P_j$ are distributed according to the decomposition $S$.

Let $(C; P_1, \ldots, P_n)$ be a fiber of $f_{0,n}$.
$C$ is a tree of rational curves.  Let $C_1, \ldots, C_c$
be its irreducible components, and $T(C_i) = \{j ; P_j \in C_i\}$.
Thus the dual graph $\tau$ of $C$ has $c$ verteces and $n$ tails, and
$\{1, \ldots, n\} = \coprod T(C_i)$.
The locus $\Delta_{\tau} \subset M_{0,n}$ which parametrizes the
$n$-pointed stable curves with the dual graph $\tau$ has codimension $c - 1$.
In fact, $\Delta_{\tau} = \bigcap_{k=1}^{c-1} \Delta_{S_k}$,
where the $S_k$ run the decompositions which corresponds to
$c-1$ decompositions of $C$ into 2 connected curves
$C^{\prime (k)}$ and $C^{\prime \prime (k)}$ so that
$S'_k = \{j ; P_j \in C^{\prime (k)}\}$ and
$S''_k = \{j ; P_j \in C^{\prime \prime (k)}\}$.

Let $d_j$ ($j = 1, \ldots, n$) be rational numbers such that
$0 < d_j \le 1$ for all $j$ and $\sum_j d_j = 2$.
Let $\Cal D = \sum_j d_j\Cal P_j$.
There is a canonical section $\omega_{\Cal D,m}$ of
\linebreak
$\Cal O_{\Cal Z_{0,n}}(mp(K_{\Cal Z_{0,n}/\Cal M_{0,n}} + \Cal D))$
over $f_{0,n}^{-1}\Cal M_{0,n}^o$.
We set $\alpha(S') = \sum_{j \in S'} d_j$
for a decomposition $S = \{S', S''\}$,
$$
\Cal F_S = \cases &(\alpha(S') -1)\Cal F_{S''} \text{ if }
\alpha(S') \ge 1 \\
&(1-\alpha(S'))\Cal F_{S'} \text{ otherwise }
\endcases
$$
and $\Cal F = \sum_S \Cal F_S$.

\proclaim{Lemma 3}
Let $C$ be any fiber of $f_{0,n}$.

(1) $(K_{\Cal Z_{0,n}} + \Cal D - \Cal F) \vert_C \equiv 0$.

(2) There exists an irreducible component $C_0$ of $C$ such that
$C_0 \not\subset \text{Supp}(\Cal F)$.
\endproclaim

\demo{Proof}
(1) We shall construct a sequence of closed subschemes $C^{(k)}$ of $C$
with contraction morphisms $C \to C^{(k)}$
for $k = 0, 1, \ldots, c$ inductively as follows.
Let $C = C^{(0)}$.
Assume that $C^{(k)}$ is already constructed for some $k$.
Then pick an irreducible component $C_{i_k}$ of $C^{(k)}$ which
intersects the rest of the curve at only 1 point.
Let $C^{(k+1)}$ be the curve obtained from $C^{(k)}$ by contracting
$C_{i_k}$ to a point,
and $C^{\prime (k)}$ the total inverse image of $C_{i_k}$ by
the natural morphism $C \to C^{(k)}$.
Then the irreducible components of $C$ other than those of $C^{\prime (k)}$
make up a connected curve $C^{\prime \prime (k)}$.
Let $S_k$ be the decomposition corresponding to the decomposition
$C = C^{\prime (k)} \cup C^{\prime \prime (k)}$, and $\tilde C^{\prime (k)}$
the fiber of $\Cal F_{S_k'}$ over a general point of $\Delta_{S_k}$.
Then we have
$((K_{\Cal Z_{0,n}} + \Cal D - \Cal F) \cdot C^{\prime (k)})
= ((K_{\Cal Z_{0,n}} + \Cal D - \Cal F) \cdot \tilde C^{\prime (k)}) = 0$.
Since $((K_{\Cal Z_{0,n}} + \Cal D - \Cal F) \cdot C_{i_k})$ is the difference
of $((K_{\Cal Z_{0,n}} + \Cal D - \Cal F) \cdot C^{\prime (k)})$ and
the sum of some of the
$((K_{\Cal Z_{0,n}} + \Cal D - \Cal F) \cdot C^{\prime (k')})$ for $k' < k$,
we obtain our claim.

(2) Let $\alpha(C^{\prime (k)}) = \sum_{P_j \in C^{\prime (k)}} d_j$, and
$k_0$ the smallest integer such that $\alpha(C^{\prime (k)}) \ge 1$.
Then $C_0 = C_{i_{k_0}}$ is the one.
\hfill $\square$
\enddemo

\proclaim{Theorem 4}
There exists a $\Bbb Q$-divisor $\Cal L$ on
$\Cal M_{0,n}$ which is effective and nef and such that
$$
K_{\Cal Z_{0,n}} + \Cal D - \Cal F \sim_{\Bbb Q}
f_{0,n}^*(K_{\Cal M_{0,n}} + \Cal L).
$$
\endproclaim

\demo{Proof}
Since $R^1f_{0,n *}\Cal O_{\Cal Z_{0,n}} = 0$,
we have $K_{\Cal Z_{0,n}} + \Cal D - \Cal F \sim_{\Bbb Q}
f_{0,n}^*(K_{\Cal M_{0,n}} + \Cal L)$ for some
$\Bbb Q$-Cartier divisor $\Cal L$.
The canonical section $\omega_{D,m}$ induces a section of
$\Cal O_{\Cal M_{0,n}}(mp\Cal L)$
over $\Cal M_{0,n}^o$.  We shall prove that this section is extended to the
whole $\Cal M_{0,n}$.
We have also to show that for any smooth
curve $B$ and any morphism $\phi: B \to \Cal M_{0,n}$, we have
$\text{deg }\phi^*\Cal L \ge 0$.
Let $(f: X \to B; P_1, \ldots, P_n)$ be a family of $n$-pointed
stable curves corresponding to $\phi$, and $F$ the pull-back of $\Cal F$
on $X$.

First, we assume that $\phi(B) \cap \Cal M_{0,n}^o \ne \emptyset$.
Let $C^1, \ldots, C^t$ be singular fibers of $f$, and $Q_{\ell} = f(C^{\ell})$
for $\ell = 1, \ldots, t$.
By Lemma 3, there exists an irreducible component $C^{\ell}_0$ of $C^{\ell}$
for each $\ell$ which is not contained in $\text{Supp}(F)$.
Let $\mu: X \to Y$ be the birational morphism obtained by contracting all the
curves on the $C^{\ell}$ except the $C^{\ell}_0$.
Then the induced morphism $g: Y \to B$ is a $\Bbb P^1$-bundle, since it has
only reduced and irreducible fibers.
Let $\bar P_j = \mu_*P_j$ and $\bar D = \sum_j d_j\bar P_j$.
The canonical section $\omega_{D,m}$ induces a section
of $\Cal O_Y(mp(K_{Y/B} + \bar D))$
which vanishes over the points $Q_{\ell}$, because some of the
sections $\bar P_j$ meet over these points.
Since $g^*\phi^*\Cal L = K_{Y/B} + \bar D$, we obtain
$\text{deg }\phi^*\Cal L \ge 0$ in this case.
The effectivity of $\Cal L$ is also proved.

Next, we consider the general case.
Let $C$ be a general fiber of $f$, and $C_0$ its irreducible component given
by the Lemma 3.  We apply the previous argument to the
irreducible component $X^{\#}$ of $X$ which contains $C_0$.
Let $\rho: X \to X^{\#}$ be the natural contraction morphism over $B$.
We define $D^{\#} = \rho_*D =
\sum_k \tilde d_k P^{\#}_k$
as follows.
The $P^{\#}_k$ are either one of the $P_j$ which are contained in
$X^{\#}$ or one of the intersection loci on $X^{\#}$
with other irreducible components of $X$.
We set $d^{\#}_k = d_j$ in the former case, and
$d^{\#}_k = \sum d_j$ where the sum is taken over those $P_j$
which are mapped to $P^{\#}_k$ by $\rho$ in the latter case.
We have $d^{\#}_k \le 1$ for all $k$, because
$C_0 \not\subset \text{Supp}(F)$.

Let $C^{\# 1}, \ldots, C^{\# t}$ be singular fibers of
$f^{\#}: X^{\#} \to B$,
$Q_{\ell} = f(C^{\# \ell})$, and $C^{\ell} = f^{-1}(Q_{\ell})$.
If $\rho^{\ell}: C^{\ell} \to C^{\# \ell}$ is the natural contraction
morphism,
then $D^{\#} \vert_{C^{\# \ell}} = \rho^{\ell}_*(D \vert_{C^{\ell}})$.
Since $d^{\#}_k \le 1$ for all $k$,
we have $C^{\# \ell} \not\subset \text{Supp}(F)$, and
there exists an irreducible component
$C^{\# \ell}_0$ of $C^{\# \ell}$ which is not contained in
$\text{Supp}(F)$ as in Lemma 3.
Let $\mu: X^{\#} \to Y^{\#}$
be the birational morphism obtained by contracting all
the curves on the $C^{\# \ell}$ except the $C^{\# \ell}_0$.
The induced morphism $g^{\#}: Y^{\#} \to B$ is a $\Bbb P^1$-bundle.
Let $\bar P_k = \mu_*P^{\#}_k$, and
$\bar D = \sum_k d^{\#}_k \bar P_k$.
The canonical section $\omega_{D^{\#}, m^{\#}}$ induces a section
of $\Cal O_{\tilde Y}(m^{\#}p^{\#}(K_{Y^{\#}/B} + \bar D))$
which vanishes over the $Q_{\ell}$.
Since $g^{\# *}\phi^*\Cal L = K_{Y^{\#}/B} + \bar D$, we obtain
$\text{deg }\phi^*\Cal L \ge 0$ again.
\hfill $\square$
\enddemo

\definition{Remark 5}
$\Cal L$ is not necessarily ample.
For example, let $X'$ and $Y''$ be trivial $\Bbb P^1$-bundles
over a proper smooth curve $B$.
Let $P_1, P_2, \Gamma'$
(resp. $\bar P_3, \bar P_4, \bar \Gamma''$) be constant sections of
$X'$ (resp. $Y''$),
and $\bar P_5$ a non-constant section of $Y''$ which intersects
the other sections transversally.
Let $\mu: X'' \to Y''$ be the blow-up at these intersection points, and
denote by $P_3, P_4, P_5, \Gamma''$ the strict transforms of
$\bar P_3, \bar P_4, \bar P_5, \bar \Gamma''$, respectively.
Let $X = X' \cup X''$ be the surface obtained by gluing $\Gamma'$ and
$\Gamma''$.
Then the natural projection $f: X \to B$ is a non-trivial family of 5-pointed
stable curves of genus 0.
Let $d_j = \frac 47, \frac 47, \frac 27, \frac 27, \frac 27$ for
$j = 1,2,3,4,5$, respectively.
Then we have $\alpha(S') = \frac 87 > 1$, and
$\text{deg }L = 0$ in this case.

But $\Cal L$ might be $\Bbb Q$-free, i.e., $\vert m\Cal L \vert$
might be free for some positive integer $m$.
If this is the case, then we do not need to assume that $X$ is affine at the
final part of Theorem 1.
\enddefinition

\proclaim{Theorem 6}
Let $f: X \to B$ be a proper surjective morphism of smooth varieties and
$P = \sum P_j$ (resp. $Q = \sum_{\ell} Q_{\ell}$) a
normal crossing divisor on $X$ (resp. $B$) such that
$f^{-1}(Q) \subset P$ and
$f$ is smooth over $B \setminus Q$ with fibers isomorphic to $\Bbb P^1$.
Let $D = \sum_j d_jP_j$ be a $\Bbb Q$-divisor on $X$,
where $d_j = 0$ is allowed,
which satisfies the following conditions:

(i) $D = D^h + D^v$ such that
$f: \text{Supp}(D^h) \to B$ is generically finite and etale
over $B \setminus Q$, and
$f(\text{Supp}(D^v)) \subset Q$.  An irreducuble
component of $D^h$ (resp. $D^v$) is called
{\it horizontal} (resp. {\it vertical}).

(ii) $0 < d_j \le 1$ if $P_j$ is horizontal and
$d_j < 1$ if vertical.

(iii) $K_X + D \sim_{\Bbb Q} f^*(K_B + L)$ for some
$\Bbb Q$-divisor $L$ on B.

Assume that there exist a finite surjective morphism
$\pi_0: \tilde B \to B$ from a smooth variety $\tilde B$ with a normal
crossing divisor $\tilde Q = \pi_0^{-1}(Q) =
\sum_m \tilde Q_m$,
a morphism $\phi_0: \tilde B \to \Cal M_{0,n}$,
and a common desingularization $\tilde X$ of $X \times_B \tilde B$ and
$\Cal Z_{0,n} \times_{\Cal M_{0,n}} \tilde B$ over $\tilde B$
as in the following commutative diagram:
$$
\CD
X      @<{\pi}<<     \tilde X     @>{\phi}>>     \Cal Z_{0,n} \\
@VfVV               @V{\tilde f}VV               @V{f_{0,n}}VV    \\
B      @<{\pi_0}<<   \tilde B     @>{\phi_0}>>     \Cal M_{0,n}
\endCD
$$
such that $\tilde P = \pi^{-1}(P) = \sum_k \tilde P_k$
is a normal crossing divisor, and
the horizontal components of $\pi^*D$ and $\phi^*\Cal D$ coincide.
Let
$$
\align
f^*Q_{\ell} &= \sum_j w_{\ell j}P_j \\
\bar d_j &= \frac {d_j + w_{\ell j} - 1}{w_{\ell j}}
\text{ if } f(P_j) = Q_{\ell} \\
\delta_{\ell} &= \text{max }\{ \bar d_j; f(P_j) = Q_{\ell}\} \\
\Delta &= \sum_{\ell} \delta_{\ell}Q_{\ell} \\
L &= M + \Delta.
\endalign
$$
Then $\delta_{\ell} < 1$ for all $\ell$,
and $M$ is effective and nef.
Moreover, if $1-w_{\ell j} \le d_j$ for some $j$,
then $0 \le \delta_{\ell}$.
\endproclaim

\demo{Proof}
The only non-trivial statement is that $M$ is effective and nef.
We define a $\Bbb Q$-divisor $\tilde D$ on $\tilde X$ by
$\pi^*(K_{X/B} + D - f^*\Delta) = K_{\tilde X/\tilde B} + \tilde D$.
We write $\tilde D = \sum_k \tilde d_k \tilde P_k$.
Let us take any $\tilde P_k$ such that $\tilde f(\tilde P_k) =
\tilde Q_m$.  Let $Q_{\ell} = \pi_0(\tilde Q_m)$ and
$P_j = \pi(\tilde P_k)$.
Let $e_m$ (resp. $e_k$) be the ramification index of $\pi_0$ (resp. $\pi$)
along $\tilde Q_m$ (resp. $\tilde P_k$).
Since $e_m = e_kw_{\ell j}$, we have
$$
\tilde d_k
= e_k(d_j - \delta_{\ell} w_{\ell j}) - (e_k - 1) + (e_m - 1)
= e_k(\bar d_j - \delta_{\ell}).
$$
Therefore, the vertical component of $- \tilde D$ over the generic point
$\tilde \eta_m$ of $\tilde Q_m$ is
effective and its support does not contain the whole fiber
$\tilde f^{-1}(\eta_{\tilde \ell})$.
So we have
$K_{\tilde X/\tilde B} + \tilde D \sim_{\Bbb Q}
\phi^*(K_{\Cal Z_{0,n}/\Cal M_{0,n}} + \Cal D - \Cal F)$,
since both hand sides are pull-backs of some $\Bbb Q$-divisors on
$\tilde B$ anyway.
Hence $\pi^*f^*(L - \Delta) \sim_{\Bbb Q} \phi^*f_{0,n}^*\Cal L$, and
$M$ is effective and nef.
\hfill $\square$
\enddemo

\definition{Remark 7}
$\Delta$ is independent of the birational model of $X$ in the following sense.
If we blow-up $X$ along the intersection of $P_{j'}$ and $P_{j''}$
over $Q_{\ell}$ and $P_j$ is the exceptional divisor, then
we have $w_j = w_{j'} + w_{j''}$ and $d_j = d_{j'} + d_{j''} - 1$.
So $\delta_{\ell}$ does not change.
\enddefinition

\definition{Example 8}
Let $g: Y \to B$ be a $\Bbb P^1$-bundle over a smooth curve with 4 sections
$\bar P_j$ $(j = 1, \ldots, 4)$.  Assume that $\bar P_1$ and $\bar P_2$
intersect at $P \in Y$ with multiplicity $n$ and that there are no other
intersections among the $\bar P_j$.
Let $\mu: X \to Y$ be the composition of $n$ blow-ups over $P$ which makes
the strict transforms $P_j = h_*^{-1}\bar P_j$ disjoint.
Let $f = g \circ \mu$ and $Q = g(P)$.  Then we have
$$
K_X + \frac 12 \sum_j P_j = \mu^*(K_Y + \frac 12 \sum_j \bar P_j) =
f^*(K_B + \frac n6 Q).
$$
Let $\rho: Z \to X$ be the double cover whose ramification divisor is
$\sum_j P_j$.  Then $K_Z = \rho^*(K_X + \frac 12 \sum_j P_j)$ and the
induced morphism $e: Z \to B$ gives an elliptic surface with a degenerate
fiber of type $I_{2n}$ over $Q$.  By the canonical bundle formula of Kodaira,
we have $K_Z = e^*(K_B + \frac {2n}{12}Q)$.
\enddefinition

\demo{Proof of Theorem 1}
Let $\mu: Y \to X$ be an embedded resolution of the pair $(X, D)$.
We write
$$
K_Y + E + F = \mu^*(K_X + D)
$$
where $E$ is the place of canonical singularities corresponding to $W$.
We may assume that there is a resolution of singularities $\nu: V \to W$
which factors $\mu: E \to W$ and
such that $f: E \to V$ and $F \vert_E$ satisfy the conditions of Theorem 6
in places of $f: X \to B$ and $D$, respectively, since
$$
\mu^*(K_X + D) \vert_E = (K_Y + E + F) \vert_E = K_E + F \vert_E.
$$
So there exist $\Bbb Q$-divisors $M$ and $\Delta$ on $V$
such that $K_E + F \vert_E = f^*(K_V + M + \Delta)$ where
$M$ is effective and nef and the coefficients of $\Delta$ are less than 1.
We put $M_W = f_*M$ and $D_W = f_*\Delta$.

Now we use the notation of Theorem 6.
By using Kawamata-Viehweg vanishing theorem, we
conclude that the natural injective homomorphism
$\Cal O_W \to \mu_*\Cal O_E(\ulcorner - F \urcorner)$
is surjective ([K3, Theorem 1.6]).
Therefore, if $\nu_*Q_{\ell} \ne 0$, then
there exists a $j$ such that $- d_j \le w_{\ell j} - 1$.
Hence $D_W$ is effective.

Assume that $X$ is affine.
Since $M$ is nef, we may assume that
$M - \epsilon \sum_{\ell} q_{\ell}Q_{\ell}$ is ample for $0 < \epsilon \ll 1$
and for some rational numbers $q_{\ell}$ such that
$q_{\ell} > 0$ (resp. $= 0$) if $\nu_*Q_{\ell} = 0$ (resp. $\ne 0$).
We take an effective $\Bbb Q$-divisor
$M' \sim_{\Bbb Q} M - \epsilon \sum_{\ell} q_{\ell}Q_{\ell}$ which
has smooth support and transversal to the $Q_{\ell}$, and
let $M'_W = \nu_*M'$.
Then we have
$$
\nu^*(K_W + M'_W + D_W)
= K_V + M' + \Delta + \epsilon \sum_{\ell} q_{\ell}Q_{\ell}
$$
since the pull-back of the right hand side by $f$ is equal to
a pull-back by $\mu$.
If $\epsilon$ is chosen small enough,
then the coefficients of the $Q_{\ell}$ on
the right hand side are less than 1, and
$(W, M'_W + D_W)$ is KLT.
\hfill $\square$
\enddemo

\definition{Remark 9}
One could use the moduli space of log surfaces [A] in order to extend our
results to the case of codimension 3.
The argument of this paper is similar to the one used in the
proof of the additivity of the logarithmic Kodaira dimension of algebraic
fiber spaces of open curves ([K1]).
The first result toward the additivity of the (non-log) Kodaira dimension used
some properties of the moduli space of stable curves ([V]),
but these were replaced later by some positivity results from the
Hodge theory in order to generalize to the case of higher dimensional fibers
([K2]).  The most important statement of Theorem 6 is that $M$ is nef,
and this might be generalized by the Hodge theoretic method to
the case of higher dimensional fibers, hence the higher codimensional
subadjunction.
\enddefinition

\enddocument